\documentclass{article}

\usepackage{amsmath}
\usepackage{url}
\usepackage{graphicx}

\begin{document}

\title{A Conjecture on the Collatz-Kakutani Path Length for the Mersenne Primes}
\author{Toru Ohira and Hiroshi Watanabe}
\date{}
\maketitle

\begin{abstract}
  We present here a new conjecture for the nature of the Mersenne prime numbers by connecting it with the
Collatz-Kakutani problem. 
By introducing a  natural \textit{path length} on the basis of the Collatz-Kakutani tree,
we conjecture that this path length of a Mersenne prime from the root of the Collatz-Kakutani tree is approximately proportional to the index of the Mersenne prime.
We also discuss difference 
of behaviors between Mersenne numbers and Mersenne primes.
\end{abstract}

\section{Introduction}

Prime numbers have attracted mathematically oriented minds.
They are fountains of the interesting problems which are left unresolved to the present.
As well as prime numbers, natural numbers themselves sometimes show
unexpected behaviors in spite of their simple appeal.
In this note,  we would like to present a simple property arising
from the combination of two unsolved problems in number theory; 
the Mersenne prime numbers and the Collatz-Kakutani conjecture.
Namely,  we show that there exists an approximate linear relation
between an index
of the Mersenne prime and its Collatz-Kakutani path length, both of which are defined in the next section.

\section{Mersenne primes and Collatz-Kakutani Conjecuture}

Let us start with a brief descriptions of the Mersenne primes and the Collatz-Kakutani conjecture\cite{mathlib}.
A Mersenne number is a positive integer given by
$$ M_n = 2^n -1,  $$
where $n$ is a positive integer,  called ``index''. A Mersenne prime is a Mersenne number that is prime.
It has been shown that if $M_n$ is a Mersenne prime,  $n$ must be prime. However,  
the converse is not true. In fact,  only forty-seven Mersenne primes have been found up to date with
the largest one given when
$n=43,112,609$. This is also the largest known prime number.
There are many fundamental questions left unresolved such as whether there 
are infinitely many Mersenne primes.
The behaviors of Mersenne primes are irregular and unpredictable as seen in
other types of prime numbers.
Therefore, the recent discoveries of the Mersenne primes were achieved with the help of
massively distributed computers, and the effort for finding new Mersenne primes
has been made continuously~\cite{gimps}. The Mersenne primes are known for their relationships
with the  perfect numbers, and they are also applied to create pseudorandom number
generators~\cite{mt}.

The Collatz-Kakutani conjecture is also one of the unsolved problems
in number theory.
This conjecture is a halting problem that the following operations will stop
for an arbitrary positive number. Consider a positive integer $X$.

If $X$ is odd,  $3X+1$ is the next integer.

If $X$ is even,  $X/2$ is the next integer.

\noindent
If we repeat this process,  it will eventually reach the halting state $X=1$
for all
positive integers. For example,  if we start from $X=7$,  
we have a sequence as 

$$7,  22,  11,  34,  17,  52,  26,  13, 40,  20,  10,  5,  16,  8,  4,  2,  1.$$

This process can be visualized as the Collatz-Kakutani tree as shown in Fig.~1A. The conjecture 
states that this tree covers all the positive integers. Even though this conjecture is unsolved,  it is generally believed as true from arguments through probability theories and
through computational verifications.

Let us introduce the ``Collatz-Kakutani path length'' $D(X)$ for a number of steps for $X$ to
reach $X=1$ through the above operations,  i.e.,  a number of operational steps needed to reach $1$ (``the root'') on the Collatz-Kakutani tree.
For example,  $D(7) = 16$. The relationship between $X$ and $D(X)$ is highly irregular and
no simple law for the path length is found\footnote{In fact, it may not
be defined for some $X$ if the Collatz-Kakutani conjecture is false}. As a result,  the plot of $D(X)$ versus
$X$ produces a irregular graph as shown in Fig.~1B.
While the general relationship between $X$ and $D(X)$ is unknown, 
there are some trivial relations such as $D(2^n) = n$.
Note that,  the number $2^n$ is just one larger than the Mersenne
number $M_n =2^n - 1$.

\begin{figure}[t]
\begin{center}
\includegraphics[width=0.95\columnwidth]{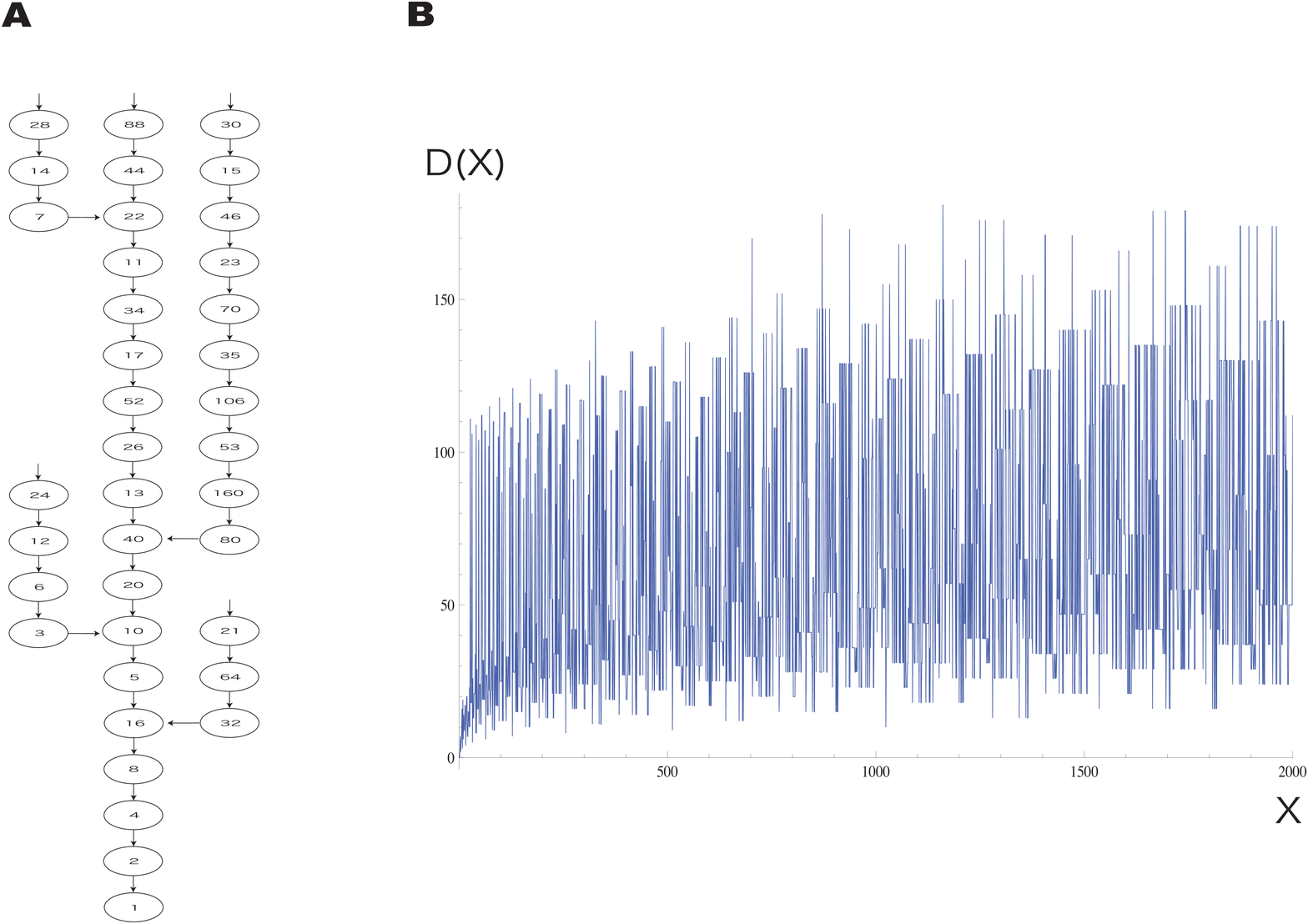}
\end{center}
\caption{A. Schematic View of Collatz-Kakutani tree. B. Plot of the path length $D(X)$ on the Collatz-Kakutani tree to the root of tree
$X=1$ for integers up to X=2000.}
\label{tree}
\end{figure}

\section{Main Results}

Our main finding to report is the fact that a path length of a Mersenne prime
is approximately proportional to its index for large n, namely, 
\begin{equation}
D(M_n) \approx 13.45 n.
\end{equation}
This is shown in Fig.~2,  which we computed up to the largest known 47th Mersenne prime,  ${M_p}(47) = M_{43112609}$. 
The behavior of the path length $D$ is non--monotonic for small indices, e.g., $D(M_{89}) > D(M_{107})$. 
We expect,  however,  
beyond $n=107$,  the path length increases monotonically with $n$ for the Mersenne primes. Also,  there seems no regular path to reach the root of the tree,  $X=1$. For example,  we have examined the path for two nearby Mersenne primes, ${M_p}(16) = M_{2203}$ and ${M_p}(17) = M_{2281}$.
While they both have the path length $D$ of approximately $30,000$,
they share the last 16 steps, namely, 

$22, 11, 34, 17, 52, 26, 13, 40, 20, 10, 5, 16, 8, 4, 2,$ and $1$.

The relations between indices of the Mersenne primes and their path lengths are summarized in Table~\ref{tbl_pathlength}.
\vspace{2em}

\begin{table}[htbp]
\footnotesize
\begin{center}
\begin{tabular}{|l|l|l|l|}
\hline
\multicolumn{1}{|l|}{Number} & \multicolumn{1}{|l|}{Indices} & \multicolumn{1}{|l|}{$D(X)$}  & \multicolumn{1}{|l|}{${D(X)/n}$}\\
\hline
1 & 2 & 7 & 3.5 \\
\hline
2 & 3 & 16 & 5.33333   \\
\hline
3 & 5 & 106 &  21.2 \\
\hline
4 & 7 & 46 &  6.57143 \\
\hline
5 & 13 & 158 & 12.1538  \\
\hline
6 & 17 & 224 & 13.1765  \\
\hline
7 & 19 & 177 & 9.31579  \\
\hline
8 & 31 & 450 &  14.5161 \\
\hline
9 & 61 & 860 &  14.0984 \\
\hline
10 & 89 & 1454 & 16.3371  \\
\hline
11 & 107 & 1441 & 13.4673  \\
\hline
12 & 127 & 1660 &  13.0709 \\
\hline
13 & 521 & 6769 &  12.9923 \\
\hline
14 & 607 & 8494 & 13.9934  \\
\hline
15 & 1279 & 17094 & 13.3651  \\
\hline
16 & 2203 & 29821 & 13.5365  \\
\hline
17 & 2281 & 30734 & 13.4739  \\
\hline
18 & 3217 & 43478 & 13.5151  \\
\hline
19 & 4253 & 55906 & 13.1451  \\
\hline
20 & 4423 & 60716 & 13.7273  \\
\hline
21 & 9689  & 129608 & 13.3768  \\
\hline
22 & 9941 &  134345 & 13.5142  \\
\hline
23 & 11213 & 153505 & 13.6899  \\
\hline
24 & 19937 & 265860 &  13.335 \\
\hline
25 & 21701 & 293161 & 13.5091  \\
\hline
26 & 23209 & 312164 & 13.4501  \\
\hline
27 & 44497 & 598067 &  13.4406 \\
\hline
28 & 86243 & 1158876 & 13.4373  \\
\hline
29 & 110503 & 1482529 & 13.4162  \\
\hline
30 & 132049 & 1771117 & 13.4126  \\
\hline
31 & 216091 & 2906179 &  13.4489 \\
\hline
32 & 756839 & 10197081 &  13.4732 \\
\hline
33 & 859433 & 11568589 & 13.4607  \\
\hline
34 & 1257787 & 16927967 & 13.4585  \\
\hline
35 & 1398269 & 18807193 & 13.4503  \\
\hline
36 & 2976221 & 40055567 & 13.4585 \\
\hline
37 & 3021377 & 40663017 & 13.4584  \\
\hline
38 & 6972593 & 93778449 & 13.4496  \\
\hline
39 & 13466917 & 181209792 & 13.4559  \\
\hline
40 & 20996011 & 282515044 & 13.4557  \\
\hline
41 & 24036583 & 323346876 & 13.4523  \\
\hline
42 & 25964951 & 349304386 & 13.4529  \\
\hline
43 & 30402457 & 409093991 & 13.456  \\
\hline
44 & 32582657 & 438465334 & 13.457  \\
\hline
45 & 371566673 & 499902411 & 13.4539  \\
\hline
46 & 42643801 & 573966881 & 13.4596  \\
\hline
47 & 43112609 & 580260946 & 13.4592  \\
\hline
\end{tabular} 
\end{center}
\caption{
Path lengths of the Mersenne primes.
}
\label{tbl_pathlength}
\end{table}
\clearpage
This linear behavior can be understood by the following heuristic arguments.
Suppose $N$ is a large random number. The standard heuristic arguments
for the path length for $N$ gives
\begin{equation}
D(N) \approx (3/(\ln 4/3)) \ln N.
\end{equation}
A Mersenne number $2^n - 1$ becomes $3 \times {2^{n-1}} - 1$ after two steps in
 Collatz-Katutani tree. Therefore, A Mersenne Number $2^n - 1$ becomes $3^n-1$ after $2n$ steps.
If we consider the number $N = 3^n-1 \sim 3^n$ to be a large random number for large $n$, then 
the path length of the number is given by
\begin{equation}
D(3^n-1) \approx (3/(\ln 4/3)) n \ln 3.
\end{equation}
Finally, we obtain the heuristic estimation of the path length of the Mersenne number to be
\begin{equation}
D(2^n -1) \approx 2n +  (3/(\ln 4/3)) n \ln 3 \approx 13.45652 n,
\end{equation}
which is consistent with our numerical results~\cite{referee1,lagarias,kontorovich}.
This arguments implies that our finding is true not only for the Mersenne primes, but also for the Mersenne numbers.
Our main conjecture, however, is that this linear relationship for the Mersenne primes 
are better than that for the general Mersenne numbers, or other general sequence of integers.

 In order to show the difference between the Mersenne numbers in general
and the Mersenne primes, the ratio of the
path length to the index, $D(M_n)/n$, 
is plotted against $n$ in Fig.~2B.  
If the distance increases linearly to the index, 
then $D(M_n)/n$
becomes a constant. Figure~2B shows that the Mersenne primes show
better linearity than the general Mersenne numbers.
\begin{figure}[h]
\begin{center}
\includegraphics[width=1.0\columnwidth]{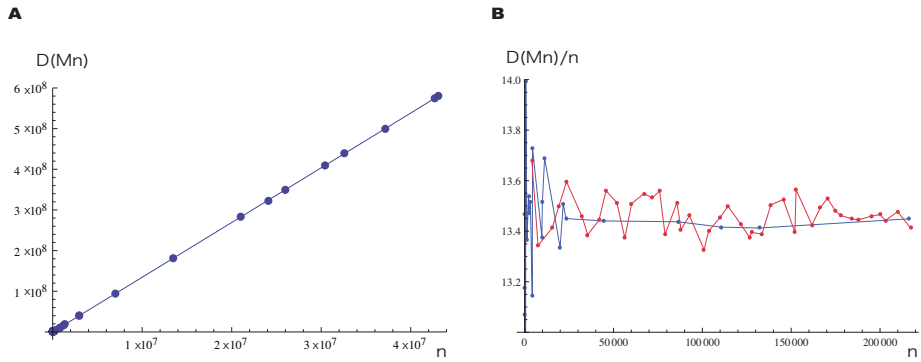}
\end{center}
\caption{A. Plot of $D(M_n)$ for Mersenne primes with index $n$ up to around ${M_p(47)}$. B. Plot of $D(M_n)/n$ 
for the Mersenne primes (Blue) and the Mersenne numbers with randomly chosen index (Red) up to around ${M_p(31)}$.}
\label{linear}
\end{figure}
\vspace{2em}

To be more quantitative,  we have performed the following statistical analyses.
We have chosen thirteen Mersenne primes from the $26^{\text{th}}$ (${M_p}(26) = M_{23, 209}$) to the $38^{\text{th}}$ (${M_p}(38) = M_{6, 972, 593}$), 
and computed the mean and the variance of $D(M_n)/n$. Let us call the prime indices of the Mersenne primes as the
Mersenne prime indices.
We compare the result with that
from the following set of thirteen points approximately in the same interval. 

A) Thirteen prime indices which is the next smallest to thirteen Mersenne prime indices.
We note the Mersenne numbers associated with this set is not prime numbers themselves.

B) Thirteen indices which are near the mid-point of two successive Mersenne prime indices.

C) Thirteen indices which are twice the Mersenne prime indices. These are all even 
numbers. We have selected some which are close to the Mersenne prime indices.

D) Thirteen indices which is on the least squre best-fit line based on the heuristics that,
given the $k$-th Mersenne prime ${M_p}(k)$,  the plot of $\log_{2} (\log_{2} ( {M_p}(k) )$ versus $k$ 
lies approximately on the straight line (Figure 3). We have selected some which are close to the Mersenne prime indices~\cite{schroeder}.
\begin{figure}[h]
\begin{center}
\includegraphics[width=0.5\columnwidth]{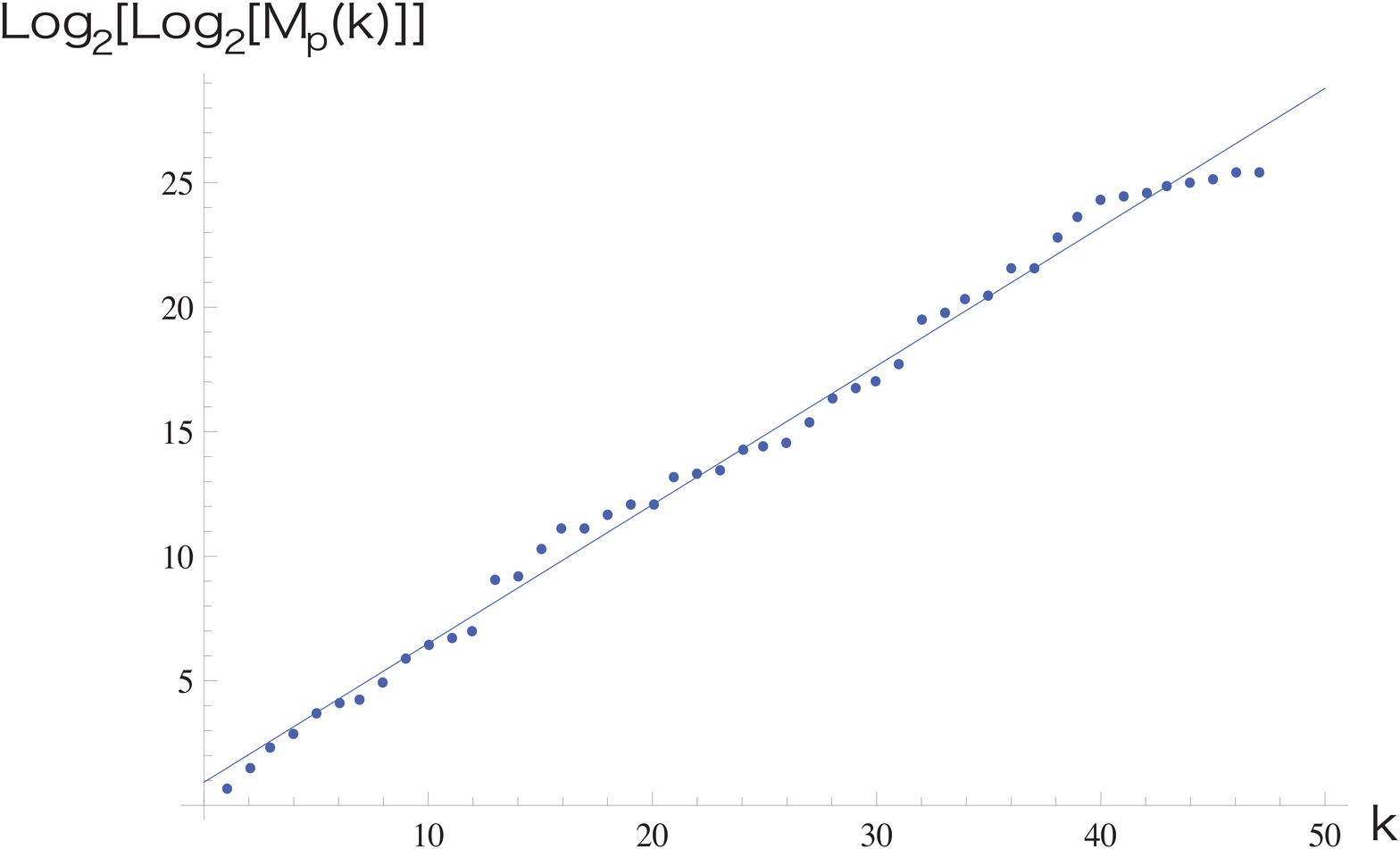}
\end{center}
\caption{Plot of $\log_{2} ( \log_{2} ( {M_p}(k) )$ versus $k$ for all the known Mersenne prime up to ${M_p(47)}$. The 
straight line is obtained by the least square error fit and $\sim 0.92757 + 0.55715 k $.}
\label{shroeder}
\end{figure}

The results of these comparison are summarized in Table~\ref{tbl_statistics}, showing that the path lengths of the Mersenne primes have the 
smaller variance than other sets.
However, this difference is not enough to single out a Mersenne prime index.

In spite of the above,  we look more detail into this property. We have examined the plot of $D(X=2^n-1)$ versus $n$
more closely and noted that they are composed of a collection of ``flat'' regions with jumps like a staircase. 
Correspondingly, the plot of ${D(X)/n}$ versus $n$ becomes collection of stripes,  which interestingly appears
near parallel. They are plotted in 
Figure 4. These plots are made for the sampled prime indices near three Mersenne prime indices,  $23209,  110503, 216091$ which give $26^{\text{th}}$, $29^{\text{th}}$, and $31^{\text{th}}$ Mersenne primes respectively. It is also interesting that there 
are occasional ``overlaps'' of ``stair steps''.


\begin{table}
\begin{tabular}{|p{16mm}|p{31mm}|p{34mm}|c|c|}
\hline
\multicolumn{1}{|c|}{Set} & \multicolumn{1}{|c|}{Indices} & \multicolumn{1}{|c|}{$D(X)$}  & \multicolumn{1}{|c|}{E$[D(X)/n]$}  & \multicolumn{1}{|c|}{Var[$D(X)/n$]}  \\
\hline
Mersenne Prime & 
23209,  44497,  86243,  110503,  132049,  216091,  756839,  859433,  1257787,  1398269,  2976221,  3021377,  6972593 & 
312164, 598067,  1158876,  1482529,  1771117,  906179,  10197081,  11568589,  16927967,  18807193,  40055567, 40663017, 93778449 & 
13.4473 & 0.0002977 \\
\hline
A & 
23227, 44501, 86249, 110527, 132059, 216103, 756853, 859447, 1257827, 1398281, 2976229, 3021407, 6972607 &
312182, 598071, 1158882, 1482553, 1771127, 2906191, 10197095, 11568603, 16928007, 18807205, 40055575, 40663047, 93778463 &
13.4460 & 0.0003194 \\
\hline
B &
22455, 33853, 65370, 98373, 121276, 174070, 486465, 808136, 1058610, 1328028, 2187245, 2998799, 4996985 &
299801, 457438, 875438, 1327329, 1633743, 2344640, 6524449, 10868120, 14246657, 17876449, 29428265, 40364153, 67195624 &
13.4485 & 0.0017853 \\
\hline
C &
22426, 43402, 46418, 88994, 172486, 221006, 264098, 432182, 1513678, 1718866, 2515574, 2796538, 5952442 &
299772, 584422, 627877, 1201650, 2320161, 2974984, 3556035, 5828307, 20384499, 23124964, 33827530, 37632788, 80085173 &
13.4618 & 0.00132591 \\
\hline
D &
20160, 43644, 64216, 94484, 139021, 204550, 442830, 651562, 958682, 1410567, 2075452, 3053739, 4493150 &
270868, 587280, 866299, 1265873, 1871202, 2748585, 5947053, 8769774, 12911249, 18991590, 27939124, 41095221, 60441877 &
13.4515 & 0.000502943 \\
\hline
\end{tabular} 
\caption{
Statistical analyses on path lengths of Mersenne primes and Mersenne numbers.
E[$x$] denotes the arithmetic average of $x$ and Var[$x$] denotes the variance of $x$, respectively.
}
\label{tbl_statistics}
\end{table}
\clearpage
\begin{figure}[h]
\begin{center}
\includegraphics[width=1.0\columnwidth]{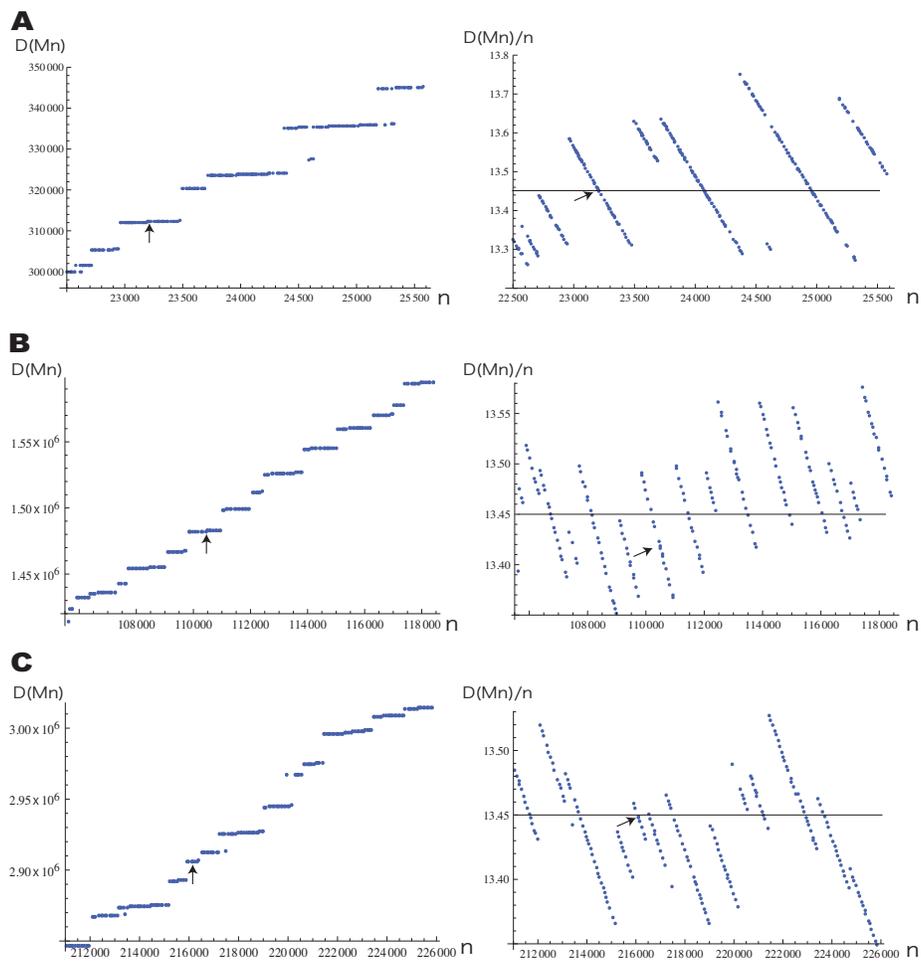}
\end{center}
\caption{Plot of $D(M_n)$ (top) and $D(M_n)/n$ (bottom) versus $n$. The Mersenne primes are 
marked by arrows, and horizontal line is at $D(M_n)/n = 13.45$. 
A. Prime indices near $n = 23209$ (${M_p(26)})$, B. Every fifth prime indices near $n = 110503$ (${M_p(29)}$),
C. Every fifth prime indices near $n = 216091$ (${M_p(31)})$.}
\label{fig4}
\end{figure}
\clearpage

The nature and mechanism of these properties appeared in Figure 4 need to be studied further.
We, nevertheless, conjecture that these qualitative behavior will continue to larger Mersenne numbers,  and the Mersenne prime
appears with its ${D(X)/n} \approx 13.45$. Unfortunately, again, this is clearly not sufficient to select out the Mersenne primes.
If this conjecture is true,  however,  we can heuristically rule out points that are further away from the horizontal line of ${D(X)/n} \approx 13.45$ from candidates of the Mersenne primes.

\section{Discussion}

We have used the Collatz-Kakutani tree as a provider of ``length'' for natural numbers, 
and obtained a conjecture for the Mersenne primes.
This may be yet another example where prime numbers can give rise to an emergence of unexpected order,  such as Ulam Spiral~\cite{ulam64, ulam67}.

Though the reason behind such simple behavior is left unclear, 
this approach of using the Collatz-Kakutani tree can be extended to other types of set of numbers known for irregularities,  possibly leading to interesting insights\cite{sinyor}.

\paragraph{Acknowledgments}  Authors would like to thank Dr.~Masaru Suzuki,  Dr.~Atsushi Kamimura, and Shigenori Matsumoto
for useful discussions. T.O. would like to thank Emeritus Prof. Yuji Ito of Keio University, for his account on his advisor, late Prof. Shizuo Kakutani of Yale University.

\bigskip

\noindent\textit{Graduate School of Mathematics, Nagoya University,  Nagoya,  Japan.\\
ohira@math.nagoya-u.ac.jp}

\bigskip

\noindent\textit{The Institute of Solid State Physics,  The University of Tokyo,  Kashiwa,  Japan.}


\begin{thebibliography}{9}
\bibitem{referee1}
This heuristic arguments of this paragraph is provided by one of the anonymous reviewers.
\end{thebibliography}
\end{document}